\documentclass[reqno]{amsart}

\usepackage{amsmath,amsthm,amssymb,amscd}
\usepackage{amsfonts}
\usepackage{rotating}
\usepackage{euscript}
\usepackage{pst-node}
\usepackage{epsfig,verbatim}

\def\be{\begin{equation}}
\def\ee{\end{equation}}

\def\C{{\mathbb C}} 
\def\f{\EuScript}
 
\def\P{{\mathbb P}}
\def\Z{{\mathbb Z}}

\def\phi{{\varphi}}

\def\deg{{\rm deg\,}}

\def\bp{\begin{proposition}}
\def\ep{\end{proposition}}

\def\bt{\begin{theorem}}
\def\et{\end{theorem}}
\def\br{\begin{remark}}
\def\er{\end{remark}}
\def\be{\begin{equation}}
\def\bee{\begin{equation*}}
\def\l{\label}

\def\ee{\end{equation}}
\def\eee{\end{equation*}}
\def\bl{\begin{lemma}}
\def\el{\end{lemma}}
\def\bc{\begin{corollary}}
\def\ec{\end{corollary}}
\def\pr{\noindent{\it Proof. }}

\def\bd{\begin{definition}}
\def\ed{\end{definition}}
\def\t{\widetilde}
\def\hat{\widehat}

\newtheorem{theorem}{Theorem}[section]
\newtheorem{lemma}{Lemma}[section]
\newtheorem{definition}{Definition}[section]
\newtheorem{corollary}{Corollary}[section]
\newtheorem{proposition}{Proposition}[section]
\newtheorem{remark}{Remark}[section]

\begin{document}
\title{Recomposing rational functions}
\author{Fedor Pakovich}

\begin{abstract} Let $A$ be a rational function.
For any decomposition of $A$ into a composition of rational functions $A=U\circ V$ the 
rational function $\t A=V\circ U$ is called an elementary transformation of $A$, and rational functions $A$ and $B$ are called equivalent  if there exists 
a chain of elementary transformations between $A$ and $B$. This equivalence relation naturally appears in the complex dynamics as a part of the problem of describing of 
semiconjugate rational functions. In this paper
we show that for a rational function $A$  its equivalence class $[A]$ 
contains infinitely many conjugacy classes if and only if 
$A$ is a flexible Latt\`es map. 
For flexible Latt\`es maps $\f L=\f L_j$
induced by the multiplication by 2 on elliptic curves 
with given $j$-invariant we provide a very precise description of  $[\f L]$.
Namely, we show that any rational function equivalent to $\f L_j$  necessarily has the form $\f L_{j'}$ for some $j'\in \C$, and that the set 
of $j'\in \C$ such that $\f L_{j'}\sim \f L_{j}$
coincides with the orbit of $j$ under the correspondence associated 
with the classical modular equation $\Phi_2(x,y)=0$. 
\end{abstract}

\maketitle

\begin{section}{Introduction}

Let $B$ be a rational function of degree at least two. The function $B$ is called semiconjugate to a rational function $A$ if the equality 
\be \l{1} A\circ X=X\circ B\ee
holds for some rational function $X$. In case if $X$ is invertible, $A$ and $B$ are called conjugate.
In terms of dynamical systems, condition \eqref{1} means that the dynamical system $A^{\circ k},$ $k\geq 1,$ on $\C\P^1$ is a factor of the dynamical system $B^{\circ k},$ $k\geq 1.$ 
The semiconjugacy is not a symmetric  relation. However, if  $B$ is semiconjugate to $A$, and  $C$ is semiconjugate to $B,$ then $C$  is semiconjugate to $A,$ since equalities \eqref{1}
and $B\circ W=W\circ C$ imply the equality  $$A\circ (X\circ W)=(X\circ W)\circ C.$$

In the recent paper \cite{semi} equation \eqref{1} was investigated at length. Roughly speaking, the main result of \cite{semi} states that  \eqref{1} holds in two cases.
In the first case, the corresponding functions $A$ and $B$ are either Latt\`es maps, or functions which can be considered as analogues of Latt\`es maps
related  to automorphism groups of $\C\P^1$  instead of automorphism groups of $\C$.  In the second case, the functions $A$ and $B$ do not possess any special properties, 
however they are  equivalent with respect to  an equivalence relation $\sim $ on the set of rational functions defined as follows.   
For any decomposition $A=U\circ V$, where $U$ and $V$ are rational functions, the 
rational function $\t A=V\circ U$ is called an {\it elementary transformation} of $A$, and rational functions $A$ and $B$ are called {\it equivalent}  if there exists 
a chain of elementary transformations between $A$ and $B$. For a rational function $A$ we will denote its equivalence class by $[A].$

The  connection between the relation $\sim$ and semiconjugacy is straightforward. Namely, for $\t A$ and $A$ as above we have:  
$$\t A\circ V=V\circ A, \ \ \ \ \ {\rm and} \ \ \ \ \ A\circ U=U\circ \t A,$$ implying inductively that 
whenever $A\sim B$ there exists $X$ such that \eqref{1} holds, and there exists $Y$ such that 
$$B\circ Y=Y\circ A.$$  Thus, if
$A\sim B$, each of  the dynamical systems $A^{\circ k},$ $k\geq 1,$ and  $B^{\circ k},$ $k\geq 1,$ 
is a factor of the other one, meaning that these systems have ``similar'' 
dynamics.
Furthermore, since for any invertible rational function $W$ the equality
$$A=(A\circ W)\circ W^{-1}$$ holds, each equivalence class $[A]$ is a union of conjugacy classes. 
Thus, the relation $\sim$ can be considered as a 
weaker form of the classical conjugacy relation.  

The main result of this paper is the following statement.

\bt \l{t1} Let $F$  be a rational function. Then its  equivalence class $[F]$ 
contains infinitely many conjugacy classes if and only if 
$F$ is a flexible Latt\`es map. 
\et

Simplest examples of  flexible
Latt\`es maps are rational functions $\f L$ 
induced by the multiplication by 2 on elliptic curves. Such a function can be defined by  
the equality 
\be \l{uxx} \wp(2z)=\f L\circ \wp(z),\ee
where $\wp(z)$ is the Weierstrass function associated with some    
lattice  $M$  of rank two in $\C$. Two such functions corresponding to lattices $M$ and $M'$ are conjugate if and only if
the elliptic curves $\C/M$ and $\C/M'$ are isomorphic. So, abusing the notation, we will denote by 
$\f L_{j}$ any Latt\`es map  
induced by the multiplication by 2 on an elliptic curve 
with given $j$-invariant.

In order to describe conjugacy classes in $[\f L_{j}]$  it is convenient to use the notion of correspondence $\f F$ associated with an affine algebraic curve $F(x,y)=0$. 
By definition, for $x_0\in \C$ an image of $x_0$  under $\f F$ is any point $y_0\in \C$ such that $F(x_0,y_0)=0$. 
More generally, $y_0\in \C$ is an image of $x_0\in \C$ under the $k$th iteration of $\f F$ if there exists a sequence $x_0,x_1,\dots , x_k=y_0$ such that $(x_{i-1},x_i)$, $i=1, \dots k,$ is a point on $F(x,y)=0.$
Considering the totalities of all images and preimages of a point $x_0$
we can define its forward,  backwards, and  full orbit  
under $\f F.$  
If $F(x,y)$ is symmetric, that is $F(x,y)=F(y,x)$, all these  orbits coincide, so we can use simply the term orbit. 

In the above notation, our main result about $[\f L_j]$ is following.

\bt \l{t2}
Any rational function equivalent to $\f L_{j}$,  $j\in \C$, has the form $\f L_{j'}$, $j'\in \C$.
Furthermore, the set 
of $j'\in \C$ such that $\f L_{j'}\sim \f L$
coincides with the orbit of $j$ under the correspondence associated 
with the classical modular equation $\Phi_2(x,y)=0$.

\et

Notice that although the expression for  $\Phi_2(x,y)=0$ is quite bulky it has a very simple parametrization by rational functions which goes back to Klein (\cite{klei}), implying that  
$\f L_{j'}\sim \f L_{j}$ if and only if  $j$ and $j'$ are in the same 
orbit of the multivalued function \be \l{kaban} \f F=\beta\circ \frac{1}{z}\circ \beta^{-1},\ee 
where $\beta$ is a rational function of degree three,  
$$\beta(z)=64\,{\frac { \left( z+4 \right) ^{3}}{{z}^{2}}}.$$

The paper has the following structure. In the second section we show that  the condition $A\sim B$ implies that $A$ and $B$ are isospectral, and 
deduce the ``only if'' part of Theorem \ref{t1} from  
the fundamental result of McMullen (\cite{Mc}) about isospectral rational functions.
In the third section we relate functional decomposition of flexible Latt\`es maps with isogenies between elliptic curves, and prove the ``if'' part of Theorem \ref{t1}.
Finally, in the fourth section we describe explicitly all functional decompositions of $\f L_j$ and prove Theorem \ref{t2}.

\end{section}

\begin{section}{Equivalence and isospectrality} 
Let $F$ be a rational function. 
Recall that two decompositions of $F$ into compositions of rational functions
$F=U\circ V$ and $F=U'\circ V'$  are called {\it equivalent} if 
there exists a M\"obius transformation $\mu$ such that 
$$U'=U\circ \mu, \ \ \ V'=\mu^{-1}\circ V.$$ Clearly, elementary transformations corresponding to 
equivalent decompositions are conjugate.  
Since 
equivalence classes of decompositions of  $F$ are in a one-to-one correspondence with 
imprimitivity systems of the monodromy group $G_F$ of $F$, this implies in particular that 
the number of conjugacy classes of rational functions obtained from $F$ by an elementary transformation is finite, and that the number of conjugacy classes in $[F]$ 
is at most countable.

Recall that a rational function $A$ is called {\it a flexible Latt\`es map} if there 
exist an elliptic curve
$\f C$ and morphisms $\phi: \f C\rightarrow \f C$ and $\pi:\f C\rightarrow \C\P^1$
such that the diagram 
\be \l{xxuuii}
\begin{CD}
\f C @>\phi>> \f C \\
@VV\pi V @VV\pi V\\ 
\C\P^1 @>A >> \ \ \C \P^1\, ,
\end{CD}
\ee
commutes, $\pi$ has degree two and satisfies $\pi(z)=\pi(-z)$, and $\phi=n z+\beta$, where $n\in \Z$ and $\beta\in \f C$ (see \cite{sildyn}, Section 6.5 and \cite{mil2}).   
In fact, 
$\beta$ necessarily satisfies the condition $2\beta=0$ on $\f C$.
Moreover, changing  $\pi(z)$ to $\pi'(z)=\pi(z+\beta)$, we see that 
the condition $\pi'(z)=\pi'(-z)$ still holds, while  \eqref{xxuuii} holds for $\phi'=n z+\beta',$ where $\beta'=n\beta.$ Thus,  if $n$ is even, we  may assume that $\beta=0$.   The complex structure of $\f C$ is completely defined by the conjugacy class of $A$, that is if $A',\f C',\pi',\phi'$ is another collection as above and $A$ is conjugate to $A'$, then $\f C$ is isomorphic to $\f C'$
(see e.g. \cite{sildyn}, Theorem 6.46).
Abusing the notation, we 
will denote by 
$A_{j}$  any Latt\`es map satisfying \eqref{xxuuii} for 
$\f C$ with given 
$j$-invariant. 

Assuming that $\f C$ is written in the Weierstrass form 
\be \l{pez} \f C:\, y^2=x^3+ax+b,\ \ \ \ \ \ \ a,b\in \C,\ee a prototypical example of a Latt\`es map is obtained for $\phi=2z$ and $\pi(x,y)=x.$ In this case, 
\be \l{wei} A=\frac{z^4-2az^2-8bz+a^2}{4z^3+4az+4b}.\ee
Notice that if $M$ is a lattice of rank two in $\C$ such that $\f C=\C/M$ and $\wp(z)$ is the corresponding Weierstrass function, then the function $A$  
is defined by the condition $\wp(2z)=A\circ \wp(z).$ 
For $\f C$ given by \eqref{pez}, $A=A_{j},$ where 
$$j=1728\frac{4a^3}{4a^3+27b^2}.$$
 
\vskip 0.2cm

Let $F$ be  a rational function of degree $d$. By definition, the {\it multiplier spectrum} of $F$ is a function which assigns to each $s\geq 1$ the unordered list of multipliers at all $d^s+1$ fixed points of $F^{\circ s}$ taken with appropriate multiplicity. 
Two rational functions  are called {\it isospectral} if they have the same   multiplier spectrum.
For example, all the functions from family \eqref{wei} have the same   multiplier spectrum (see e.g. \cite{sildyn}, Example 6.49). Nevertheless, by  the following result of McMullen 
such a  situation is exceptional (see \cite{Mc}, \cite{mil2}, \cite{sildyn}).

\bt{\rm (McMullen)} The conjugacy class of any rational function $F$  which is not a flexible Latt\`es map is defined up to finitely many choices by its
multiplier spectrum. \qed
\et

The equivalence $\sim$ and the isospectrality are closely related as the following lemma shows.

\bl  \l{uv} Let $U$ and $V$ be rational functions. Then the rational functions $U\circ V$ and $V\circ U$ are isospectral.  
\el
\pr 
Since the equality $$(U\circ V)^{\circ l}(z_0)=z_0$$  implies the equality 
$$(V\circ U)^{\circ l}(z_1)=z_1,$$ where $z_1=V(z_0),$ the function  $V$ maps periodic points of $U\circ V$ to  periodic points of $V\circ U$. Furthermore, the period of $V(z_0)$ divides the period of $z_0$. Similarly, 
$U$ maps periodic points of $V\circ U$  to  periodic points of $U\circ V$. Since the composition $U\circ V$   maps bijectively periodic points of $U\circ V$ of period $l$ to themselves, this implies 
 that  $V$  maps bijectively periodic points of $U\circ V$ of period $l$ to  periodic points of $V\circ U$ of period $l$.

Finally, since by the chain rule 
$$((U\circ V)^{\circ l})^{\prime}(z_0)=((U\circ V)^{\circ l-1}\circ U)^{\prime}(z_1)\circ V^{\prime}(z_0)$$
and 
$$((V\circ U)^{\circ l})^{\prime}(z_1)= V^{\prime}((U\circ V)^{\circ l-1}\circ U)(z_1))\circ((U\circ V)^{\circ l-1}\circ U)^{\prime}(z_1),$$
it follows from $$((U\circ V)^{\circ l-1}\circ U)(z_1)=z_0$$ that 
$$ ((U\circ V)^{\circ l})^{\prime}(z_0)=((V\circ U)^{\circ l})^{\prime}(z_1).\eqno{\Box}$$

\bc \l{uv+} Let $A$ and $B$ be rational functions such that $A\sim B$. Then $A$ and $B$ are isospectral. 
\ec
\pr By definition, $A\sim B$ if $B$ is obtained from $A$ by a chain of elementary transformations.
On the other hand, any such transformation leads to an isospectral function by Lemma \ref{uv}.  \qed

\vskip 0.2cm

It is clear that the McMullen theorem combined with Corollary \ref{uv+} proves the ``only if'' part of Theorem \ref{t1}. Notice however that the number 
of conjugacy classes in an equivalence class $[F]$ can be arbitrary large (see \cite{semi}). 
The 
proof of the ``if'' part of Theorem \ref{t1} is given in the next section. 
\vskip 0.2cm

Notice that  isospectral $A$ and $B$ are not necessary   equivalent. Say, all functions \eqref{wei} cannot be equivalent since any equivalence class contains at most countably many conjugacy classes. Nevertheless, to our best knowledge all known examples of isospectral rational functions are obtained either from flexible Latt\`es maps, or from rigid Latt\`es maps (see \cite{Mc}, \cite{mil2}, \cite{sildyn}), or else from elementary transformations. Thus, it is still possible that if $A$ and $B$ are not Latt\`es maps, then the fact that $A$ and $B$ are isospectral implies that $A\sim B$. A comprehensive description of relations  between
the isospectrality and the equi\-valence $\sim$ 
seems to be a very interesting problem.

 Notice also that if $A$ is a {\it polynomial}, then the finiteness of $[A]$ can be established without using the McMullen theorem; see Corollary 5.8 in \cite{pj},   
and also the paper \cite{ms} 
 using the notion of ``skew twist equivalence'' which is essentially coincides with the equivalence $\sim$ in the setting considered here. 
The approach of the paper \cite{ms} is based on the theory of decomposition of polynomials developed by Ritt \cite{r1}, while the method of \cite{pj} relies on the results of \cite{p1}
about polynomials sharing preimages of compact sets. However, methods of both these papers are restricted to the polynomial case only. 

\end{section}

\begin{section}{Proof of Theorem \ref{t1}}
Consider first flexible Latt\`es maps $A=A_j$ defined by the diagram 
\be \l{beg}
\begin{CD}
\f C @>n z>> \f C \\
@VV \pi  V @VV \pi V\\ 
\C\P^1 @>A >> \ \ \C \P^1\, 
\end{CD}
\ee
where $n\in \Z.$ In order to prove that $[A_j]$ 
contains infinitely many conjugacy classes we will use the relation between functional decompositions of $A_j$ and isogenies between elliptic curves.
We start from recalling some basic definitions and results concerning isogenies 
and the modular equation. Abusing the notation, below we will use the symbol $j$ in two possible meanings: for a value of the elliptic modular function $j(\tau)$
 of weight zero for $SL(2, \Z)$  on the upper half-plane $\mathbb H,$ and 
for a value of the $j$-invariant of elliptic curve $\f C=\C/L_{\tau}$,
where $L_{\tau}$ is a lattice in $\C$ generated by $1$ and $\tau\in \mathbb H.$

Let $\f C$ and $\t{\f C}$ be elliptic curves over $\C$ with corresponding neutral elements $O$ and $\t O$. An {\it isogeny} between $\f C$ and $\t{\f C}$ is a morphism $\psi:\f C\rightarrow \t{\f C}$ such that $\psi(O)=\t O$. Such a morphism is necessarily  a homomorphism of groups. A kernel  $\Gamma$  of a non-zero isogeny $\psi: \f C\rightarrow \t{\f C}$ is a subgroup of finite order in  
$\f C$, and for any  subgroup of finite order $\Gamma$ there exists a unique isogeny  $\psi:\f C\rightarrow \t{\f C}$ such that $\ker \psi=\Gamma$. For any elliptic curve $\f C$ and integer $n$ the multiplication by $n$ on $\C$ projects to an isogeny $[n]:\f C\rightarrow \f C$ of degree $n^2$ with 
kernel consisting of points whose order divides $n$.  Furthermore, for any isogeny 
$\psi:\f C\rightarrow \t{\f C}$ of degree $n$ there exists 
 a unique {\it dual} isogeny  $\widehat\psi:\, \t{\f C}\rightarrow \f C$ such that \be \l{dv} \widehat \psi \circ \psi =[n]\ee on $\f C$, 
and
\be \l{dv1}  \psi \circ\widehat \psi =[n]\ee on $\t{\f C}$ 
(for the proofs of the above facts see e.g. \cite{sil}, Chapter III). An isogeny whose kernel is a cyclic group of order $n$ is called $n$-cyclic.

Below, we mostly will consider specific $n$-cyclic isogenies
defined as follows.
Let  $\f C_1=\C/L_{\tau}$ be an elliptic curve and  $n$ an integer. Then  
the multiplication by $n$ on $\C$ projects to an $n$-cyclic isogeny 
\be \l{i1} \psi_n:\f C_{1} \rightarrow \f C_{2},\ee where  $\f C_2=\C/L_{n\tau}$. 
The dual isogeny \be \l{i2} \hat\psi_n:\f C_{2} \rightarrow \f C_{1}\ee is the projection of the identical map on $\C$.

There exists a polynomial in two variable 
\be \l{mc} \Phi_n(x,y)=0\ee with integer coefficients, called the {\it modular equation}, having the following property:
 if $\f C_1$ and $\f C_2$ are two elliptic curves with $j$-invariant $j_1$ and $j_2$, then an  $n$-cyclic isogeny  $\f C_{1}\rightarrow \f C_{2}$  exists 
if and only if $(j_1,j_2)$ is a point of curve \eqref{mc}  (see \cite{lang}, Chapter 5). In particular,
since \eqref{i1} is an $n$-cyclic isogeny between elliptic curves with $j$-invariants $j(\tau)$ and $j(n\tau),$
for any $\tau\in \mathbb H$ the equality 
\be \l{mc1} \Phi_n(j(\tau),j(n\tau))=0\ee holds. 

\vskip 0.2cm

Let $\f C$ be an elliptic curve and $A=A_j$ a Latt\`es map satisfying \eqref{beg}. Further, let  $\Gamma$ be a subgroup  of $\f C$ and  $\psi:\f C\rightarrow \t{\f C}$ an isogeny such that $\ker \psi=\Gamma$. 
Since  $\psi$  is a homomorphism, the equality $\psi(-x)=-\psi(x)$ holds,  implying that 
there exists a rational function $V_{\Gamma}$ such that 
the diagram
$$
\begin{CD}
\f C @>\psi>> \t{\f C} \\
@VV \pi  V @VV \t{\pi} V\\ 
\C\P^1 @>V_{\Gamma} >> \ \ \C \P^1\, 
\end{CD}
$$
commutes. Similarly, for a dual isogeny $\widehat\psi:\, \t{\f C}\rightarrow \f C$  there exists   
a rational function $U_{\Gamma}$ such that 
 the diagram 
$$
\begin{CD}
\t{\f C} @>\widehat \psi>> \f C \\
@VV \t\pi V @VV \pi  V\\ 
\C\P^1 @>U_{\Gamma} >> \ \ \C \P^1\, 
\end{CD}
$$ commutes. Thus, to any subgroup $\Gamma$ of $\f C$ corresponds a decomposition 
\be \l{gam} A_j=U_{\Gamma}\circ V_{\Gamma}.\ee  
In particular, 
 the equality 
$$\widehat\psi_n \circ \psi_n =[n]$$  gives rise to 
a decomposition 
 \be \l{dv2} A_j=U_{n}\circ V_{n}.\ee 
Notice that explicit expressions for $U_{\Gamma}$ and $V_{\Gamma}$ 
can be deduced from V\'elu's formulas for isogenies $\psi: \f C \rightarrow \t{\f C}$ with given $\f C$ and $\ker \psi$ 
(see \cite{velu}).

 Equality \eqref{dv1} implies that an elementary transformation of $A_j$ corresponding to decomposition \eqref{gam} also is a flexible Latt\`es map $A_{j'}$, where $j'$ is $j$-invariant of $\t{\f C}.$ Furthermore, if $\Gamma$ is a cyclic group of order $n$, then the corresponding values of $j'$ are described by the condition that  $(j,j')$ is a point of \eqref{mc}.   Clearly, 
in order to prove that $[A_j]$ has infinitely many conjugacy classes 
it is enough to prove that we can obtain infinitely many  conjugacy classes using chains of elementary transformations arising from decompositions \eqref{dv2} only. 
Moreover,  
since conjugate Latt\`es maps correspond to isomorphic elliptic curves, 
it follows from  \eqref{mc1} that it is enough to show that for any $\tau\in \mathbb H$ 
the sequence $j(n^k\tau)$ takes infinitely many values.

In order to prove the last statement recall that the Fourier expansion for $j(\tau)$ in  $q=e^{2\pi i \tau}$ is  
$$j(\tau)=\frac{1}{q}+744+196884q+\dots .$$
Further, 
$$q_k=e^{2\pi i (n^k\tau)}\to 0, \ \ \ {\rm as} \ \ \  k\to \infty,$$ 
since $\Im(\tau)>0.$ 
Therefore, 
$$j(n^k\tau)\to \infty, \ \ \ {\rm as} \ \ \   \ k\to \infty,$$ 
implying that  $j(n^k\tau)$ takes infinitely many distinct values. 
This proves  the ``if'' part of Theorem \ref{t1} for Latt\`es maps $A$ given by \eqref{beg}.

\vskip 0.2cm

Consider now flexible Latt\`es maps $\t A=\t A_j$ defined by the diagram 
\be \l{beg2}
\begin{CD}
\f C @>nz+\beta>> \f C \\
@VV \pi  V @VV \pi V\\ 
\C\P^1 @>\t A >> \ \ \C \P^1\, ,
\end{CD}
\ee
where $n\in \Z$ is odd and $\beta$ is a
point of order two on $\f C$. 
Set $\phi_{\beta}=z+\beta$.
Since  $2\beta=0$ and $n$ is odd,
\be \l{ssvv} nz+\beta=nz \circ \phi_{\beta}\ee
on $\f C$.
Further, since $\phi_{\beta}(-z)=-\phi_{\beta}(z)$ on $\f C$, there exists a M\"obius transformation $\mu$ which makes the diagram 
$$
\begin{CD}
\f C @>\phi_{\beta}>> \f C \\
@VV\pi V @VV\pi V\\ 
\C\P^1 @>\mu  >> \ \ \C \P^1
\end{CD}
$$
commutative. Therefore, the Latt\`es map $\t A=\t A_j$ 
and the Latt\`es map  $A=A_j$ defined by diagram 
\eqref{beg} are related by the equality 
\be \l{sh} \t A_j=A_j\circ \mu,\ee where $\mu$  is 
a M\"obius transformation.

Clearly, an elementary transformation $A_j\rightarrow A_{j'}$ corresponding to decomposition  \eqref{dv2} 
induces an
elementary transformation
\be \l{et}  U_n\circ V_n\circ \mu\rightarrow V_n\circ \mu\circ U_n\ee 
of $\t A_j$.
Moreover, since the composition 
$$\f C_2 \overset{\hat{\psi}_n}{\longrightarrow} \f C_1\overset{\phi_{\beta}}{\longrightarrow}  \f C_1$$
equals the composition 
$$\f C_2 \overset{\phi_{n\beta}}{\longrightarrow} \f C_2\overset{\hat{\psi}_n}{\longrightarrow}  \f C_1,$$ 
and $\phi_{n\beta}(-z)=-\phi_{n\beta}(z)$ on $\f C_2$, 
considering  
the commutative diagram 
$$
\begin{CD}
\f C_2 @>\hat\psi_n>> \f C_1 @>\phi_{\beta}>> \f C_1 @>\psi_n>> \f C_2\\
@VV\pi_2 V @VV\pi_1 V @VV\pi_1 V  @VV\pi_2 V \\ 
\C\P^1 @>U_n  >> \ \C \P^1 @>\mu  >> \ \C \P^1 @>V_n  >>\  \C \P^1
\end{CD}
$$
we see that $$(V_n\circ \mu)\circ U_n=\mu'\circ V_n\circ U_n,$$  
where  $\mu'$ is a M\"obius transformation satisfying 
$$
\begin{CD}
\f C_2 @>\phi_{n\beta}>> \f C_2 \\
@VV\pi_2 V @VV\pi_2 V\\ 
\C\P^1 @>\mu'  >> \ \ \C \P^1\, .
\end{CD}
$$
Thus, if $A_j\rightarrow A_{j'}$ is an elementary transformation  corresponding to decomposition  \eqref{dv2}, 
then elementary transformation of $\t A_j$ in the right part of \eqref{et}  is a Latt\`es map $\t A_{j'}$ defined by the diagram 
$$ \begin{CD}
\f C_2 @>nz+\beta'>> \f C_2 \\
@VV \pi_2  V @VV \pi_2 V\\ 
\C\P^1 @>\t A_{j'} >> \ \ \C \P^1\, , \end{CD}
$$ 
where $\beta'=n\beta$ is a points of order two on $\f C_2.$ 
Therefore,  a chain of  elementary transformations 
$$A_j \rightarrow A_{j_1}\rightarrow A_{j_2}\rightarrow \dots  $$ with infinitely many different $j$  arising from decompositions \eqref{dv2}
 induces a similar chain 
$$\t A_j \rightarrow \t A_{j_1}\rightarrow \t A_{j_2}\rightarrow \dots \, .$$  This finishes the proof of Theorem \ref{t1}.

\end{section}

\begin{section}{Decompositions of $\f L_j$}
In this section we provide some explicit formulas illustrating constructions from the previous section in the simplest case 
where considered Latt\`es maps are 
defined by the diagram 
\be \l{beg3}
\begin{CD}
\f C @>2z >> \f C \\
@VV \pi  V @VV \pi V\\ 
\C\P^1 @>\f L >> \ \ \C \P^1\, ,
\end{CD}
\ee
and prove Theorem \ref{t2}. 
In order to reduce the number of parameters, we will write elliptic curves in the Legendre form 
\be \l{leg} \f E_{\lambda}:\, y^2=x(x-1)(x-\lambda), \ \ \ \ \ \ \lambda \in \C\setminus\{0,1\}.\ee 
When working with the explicit expression for $\f L=\f L_j$ 
we will use the notation   
$$L_{\lambda}=\frac{1}{4}\,{\frac { \left( {z}^{2}-\lambda\right) ^{2}}{z \left( z-1 \right) 
 \left( z-\lambda\right) }}.$$ So, $L_{\lambda}=\f L_j$, where
$$j=256\frac{(\lambda^2-\lambda+1)^3}{\lambda^2(\lambda-1)^2}.$$ 
The next result describes explicitly all equivalence classes of decompositions of $L_{\lambda}$.

\bt Any decomposition of  $L_{\lambda}$ into a composition of rational functions of degree greater than one is equivalent to one of the following decompositions 

$$L_{\lambda}(z)=
\left(\frac{1}{4}\,{\frac {{z}^{2}-4\,\lambda}{z-\lambda-1}}\right) \circ 
\left(z+{\frac {\lambda}{z}}\right)
$$

\be \l{isog} L_{\lambda}(z)= 
\left(\frac{1}{4}\,{\frac {{z}^{2}+2\,z+1}{z+1-\lambda}}\right) \circ \left(z+{\frac {1-\lambda}{z-1}}\right)
\ee

$$L_{\lambda}(z)=
\left(\frac{1}{4}\,{\frac {{\lambda}^{2}+2\,\lambda z+{z}^{2}}{z+\lambda-1}}\right)\circ \left(z+{\frac {{\lambda}^{2}-\lambda}{z-\lambda}}\right)\,.
$$
These decompositions are not equivalent and have form \eqref{gam}, where $\Gamma$ runs over cyclic subgroups of order two in $\f C.$
\et 
\pr 
Recall that 
equivalence classes of decompositions of a rational function $F$ are in a one-to-one correspondence with 
imprimitiivity systems of the monodromy group $G_F$ of $F$. Namely, if $z_0$ is a non-critical value of $F$, and $G_F$ is realized as 
a permutation groups acting on the set $F^{-1}\{z_0\},$ then to an equivalence class of a decomposition  $F=U\circ V$,  corresponds 
an imprimitivity system of $G_F$ consisting of 
$d=\deg A$ blocks $V^{-1}\{t_i\},$ $1\leq i \leq d,$ where $\{t_1,t_2,\dots, t_{d}\} =U^{-1}\{z_0\}$. 

It is 
clear that for any decomposition $L_{\lambda}=U\circ V$ with $\deg U>1,$ $\deg V>1$ the equalities $\deg V=2,$ 
$\deg U=2$ hold. Therefore, if  $c$ is a non-critical value of $L_{\lambda}(z)$ and $G_{L_{\lambda}}$ is realized as 
a permutation groups acting on $L_{\lambda}^{-1}\{c\}=\{z_0,z_1,z_2,z_3\},$
then to each equivalence class of decompositions of $L_{\lambda}$  corresponds a block of size two containing the point 
$z_0$, say. Since there might be 
at most three such blocks, namely,  $\{z_0,z_1\},$  $\{z_0,z_2\},$ and  $\{z_0,z_3\}$, there exist  
at  most three non-equivalent decompositions of $L_{\lambda}$.

Prove now that the decompositions 
$$L_{\lambda}=C_i\circ D_i, \ \ \  \ i=1,2,3,$$ 
given by formulas 
\eqref{isog} are not equivalent. 
Since the set 
$$L_{\lambda}^{-1}\{\infty\}=\{\infty,0,1,\lambda\}$$ consists of four different points, the point $c=\infty$ is a non-critical value of $L_{\lambda}$ 
and we can assume that  $G_{L_{\lambda}}$ acts on the set $L_{\lambda}^{-1}\{\infty\}$. Furthermore, 
$$C_1^{-1}\{\infty\}=\{\lambda+1,\infty\}, \ \ \  C_2^{-1}\{\infty\}=\{\lambda-1,\infty\},  \ \ \  C_3^{-1}\{\infty\}=\{1-\lambda,\infty\},$$ 
implying that the blocks containing the point $z_0=\infty$  of the corresponding imprimitivity systems are
\be \l{b} \{0,\infty\}, \ \ \ \ \{1,\infty\},  \ \ \ \ \{\lambda,\infty\}.\ee
Since these blocks are different,  decompositions 
\eqref{isog} are not equivalent.

Finally, since  any element of the function field of $\f E_{\lambda}$ has the form $p(x)+q(x)y$, where $p$ and $q$ are rational functions, 
it follows from $\pi(-z)=\pi(z)$ and $\deg \pi=2$ that $\pi=\mu\circ x$ for some M\"obius transformation $\pi$. Therefore, 
without loss of generality we may assume that $\pi(z)=x$.  
In this case,  $\pi(O)=\infty$ and $\pi$ is a bijection between the set of points satisfying $2\beta =0$ on $\f C$ and the set $L_{\lambda}^{-1}\{\infty\}$.
This implies that 
the blocks containing $z_0=\infty$ corresponding to decompositions \eqref{gam}, where $\Gamma$ is a cyclic subgroup of order two,
are exactly the two-element sets listed in \eqref{b}, that is the blocks corresponding to decompositions \eqref{isog}. Therefore,  decompositions \eqref{gam} are equivalent to decompositions \eqref{isog}.
\qed

\vskip 0.2cm
Notice that the images of the isogenies corresponding to left parts of decompositions \eqref{isog} do not have Legendre form \eqref{leg}.  Therefore, elementary transformations  corresponding to decompositions \eqref{isog} are not {\it equal} to 
functions $L_{\lambda}(z)$  but only {\it conjugate} to such functions.

\bc \l{coro} 
Any elementary transformations of the function $\f L_{j},$ $j\in \C,$ has the form $\f L_{j'}$, $j'\in \C,$ where values of $j'$ are defined by the condition $\Phi_2(j,j')=0.$  
\ec

\pr Indeed, if $L_{\lambda}=U\circ V$ is a decomposition such that one of the functions $U$ and $V$ is invertible, then the corresponding elementary transformation leads to a  function conjugate to $L_{\lambda}$.
On the other hand, any   decomposition of  $L_{\lambda}$ into a composition of rational functions of degree greater than one
is equivalent to one of decompositions \eqref{isog}. 
 \qed

\vskip 0.2cm

Clearly, Corollary \ref{coro} implies Theorem \ref{t2}. Furthermore, since the modular equation  $\Phi_2(x,y)=0$ is given by the equation 
\begin{multline}\l{enot}
-x^2y^2+x^3+y^3+2^4\cdot 3\cdot 31\,xy(x+y)+3^4\cdot 5^3\cdot 4027\, xy\\ -2^4\cdot 3^4\cdot 5^3(x^2+y^2)+ 2^8\cdot 3^7\cdot 5^6(x+y)-2^{12}\cdot 3^9\cdot 5^9=0
\end{multline} 
and can be parametrized by the rational functions 
$$x=64\,{\frac { \left( j+4 \right) ^{3}}{{j}^{2}}}, \ \ \ \ \ \ \ \ y=64\,{\frac { \left( j+4 \right) ^{3}}{{j}^{2}}}\circ \frac{1}{j},$$ 
the correspondence $\f F$ associated with \eqref{enot} has form \eqref{kaban}.


\end{section}

\bibliographystyle{amsplain}

\end{document}